\documentclass{amsart}

\usepackage{amsmath,amsthm}
\usepackage{amsfonts,amssymb}

\newtheorem{thm}{Theorem}[section]

\newtheorem{lem}[thm]{Lemma}
\newtheorem{prop}[thm]{Proposition}

\theoremstyle{remark}

 \def\xb{{\mathbf x}}
 \def\yb{{\mathbf y}}

 \def\CA{{\mathcal A}}
 
 \def\CD{{\mathcal D}}

 \def\CI{{\mathcal I}}

 \def\CO{{\mathcal O}}
 
 \def\CR{{\mathcal R}}
 
 \def\CV{{\mathcal V}}

 \def\NN{{\mathbb N}}

 \def\RR{{\mathbb R}}

        \def\proj{\operatorname{proj}}
        
        \def\supp{\operatorname{supp}}

\newif\ifpdf
\ifx\pdfoutput\undefined
  \pdffalse
\else
  \pdftrue
\fi

\ifpdf
  \usepackage[pdftex]{graphicx}
  \DeclareGraphicsExtensions{.pdf,.jpg,.png}
\else
  \usepackage{graphicx}
\fi

\begin{document}

\title[Fast OPED algorithm for reconstruction of images]
{Fast OPED algorithm for reconstruction of images from Radon data}

\author{Yuan Xu}
\address{Department of Mathematics University of Oregon
    Eugene, Oregon 97403-1222.}
  \email{yuan@math.uoregon.edu}
\author{Oleg Tischenko}
\address{Institute of Radiation Protection\\
GSF - National Research Center for Environment and Health\\
D-85764 Neuherberg, Germany}
\email{oleg.tischenko@gsf.de}

\date{\today}
\keywords{OPED algorithms, reconstruction of images,  Radon data, polynomials of two variables}
\subjclass{42C15, 65D15}
\thanks{ The first author  was partially supported 
by the National Science Foundation under Grant DMS-0604056}

\begin{abstract}
A fast implementation of the OPED algorithm, a reconstruction 
algorithm for Radon data introduced recently, is proposed and 
tested. The new implementation uses FFT for discrete sine 
transform and an interpolation step. The convergence of the 
fast implementation is proved under the condition that the 
function is mildly smooth. The numerical test shows that 
the accuracy of the OPED algorithm changes little when the 
fast implementation is used.
\end{abstract}

\maketitle

\section{Introduction}
\setcounter{equation}{0}

Reconstruction images from Radon data is the central theme in 
medical imaging. The dominating reconstruction method, 
or algorithm, currently used for this task is FBP (filtered 
backprojection) method. A new reconstruction algorithm, called 
OPED,  was introduced recently as a possible alternative in 
\cite{X05} and studied further in \cite{XTC, XTC2}. OPED is based 
on orthogonal polynomial expansions on the disk, which is 
fundamentally different from FBP. Our study has indicated 
a number of potential advantages of OPED vs FBP.  The 
purpose of this paper is to show that OPED can compete with 
FBP in speed as well. 

Fast reconstruction is one of the reasons that FBP is widely used. 
The structure of the FBP algorithm contains a discrete convolution, 
which can be evaluated with FFT (fast Fourier transform). Together 
with an interpolation step, this ensures that FBP algorithm can be 
implemented efficiently. In contrast, the OPED algorithm does not 
use discrete convolutions. Nevertheless, it uses Chebyshev 
polynomials of the second kind, which turns out to allow an 
implementation with similar combination of FFT and interpolation, 
except that we will use fast discrete sine transforms instead. 
The result is a fast implementation of OPED algorithm, comparable 
to the implementation of FBP in the number of operations. 

The additional interpolation step will introduce an error in the 
reconstruction. It turns out, however, that the additional error 
is small. The convergence of the OPED algorithm is proved 
theoretically under mild condition on the function that represents 
the image. Treating the approximation provided by the algorithm
as an operator, denoted by $\CA_{2m}$, from the space of 
continuous functions to itself, the proof in \cite{X05} amounts to 
show that the operator norm of $\CA_{2m}$ is exactly 
$\|\CA_{2m}\| = \CO(m \log(m+1))$, where $2m+1$ is the 
number of views of the Radon projections.  It turns out, as we 
shall show, that the order of the norm remains to be the same 
even when the interpolation step is used. The OPED algorithm 
preserves polynomials of degree $2m-1$, so that the order of 
$\|\CA_{2m}\|$ implies that the algorithm converges uniformly 
if the function is mildly smooth, say, if the function has continuous 
second order derivatives. The OPED with interpolation step no 
longer preserves polynomials. Nevertheless, we shall show 
that the convergence still holds for smooth functions, say for 
functions that have fourth continuous derivatives. We also test 
the fast  implementation of OPED numerically. The result shows 
that the speed of fast implementation is improved in the order of
magnitude, whiles the quality of the reconstruction changes little 
even for images that have sharp jumps. 

The paper is organized as follows. We describe OPED algorithm 
and its fast implementation in the following section. The 
convergence of the fast OPED algorithm is established in 
Section 3. The numerical testing and examples are given in 
Section 4. 

\section{OPED algorithm and its fast implementation}
\setcounter{equation}{0}

\subsection{OPED algorithm} A two dimensional image on the 
unit disk $B = \{(x,y): x^2+y^2 \le 1\}$ is represented by a 
function $f(x,y)$ defined on $B$. A Radon projection of $f$ is a 
line integral,
$$
\CR f(\theta, t) :=  \int_{I(\theta,t)} f(x,y) dx dy,
\qquad 0 \le \theta \le 2\pi, \quad -1\le t \le 1,
$$
where $I(\theta,t) = \{(x,y): x \cos \theta + y \sin \theta = t\} \cap B$
is a line segment inside $B$. The essential problem of image 
reconstruction from Radon data is to recover the image 
represented by $f(x,y)$ from a finite collection of its Radon 
projections. 

OPED is a reconstruction algorithm based on orthogonal 
expansion on the disk. Let $\CV_n(B)$ denote the space of 
orthogonal polynomials of degree $n$ on $B$ with respect to 
the Lebesgue measure.  A function in $L^2(B)$ can be expanded 
in terms of orthogonal polynomials, that is, 
\begin{equation} \label{OPexpan}
 f(x) = \sum_{k=0}^\infty \proj_k f(x), \qquad 
        \proj_k: L^2(B) \mapsto \CV_n(B). 
\end{equation} 
The partial sum of this expansion is $S_n f(x) = \sum_{k=0}^n \proj_k f(x,y)$, 
which provides a natural approximation to $f$.  It is shown in \cite{X05} that 
the partial sum can be expressed directly in terms of Radon projections: 
\begin{equation} \label{S2m}
   S_{2m} f(x,y) =  \frac{1}{2m+1} \sum_{\nu =0}^{2m} 
        \frac{1}{\pi} \int_{-1}^1 \CR f(\phi_\nu, t)  \Phi_\nu(t;x,y) dt 
\end{equation}
where $\phi_\nu = \frac{2 \pi \nu}{2m+1}$,
$$
   \Phi_\nu(t;x,y) = \frac{1}{2m+1} \sum_{k=0}^{2m}(k+1)U_k(t) 
            U_k(x\cos\phi_\nu+ y\sin\phi_\nu),
$$
and $U_k(t)$ denote the Chebyshev polynomial of the second kind, 
\begin{equation} \label{Cheby} 
U_k(t) = \frac{\sin(k+1) \theta}{\sin \theta}, \qquad t = \cos \theta.
\end{equation}
It turns out that the formula \eqref{S2m} can be found implicitly in 
\cite{LS} (see (5.9), (4.3) and (3.7) there). Furthermore, a double
integral expression for $S_n f$ can be found in \cite{BG,P} from
which \eqref{S2m} can be deduced by a quadrature on the unit
circle. 

Using a quadrature to discretize the integral over $t$ gives an 
approximation to $f$ based on discrete Radon data. A particular 
choice of the quadrature formula gives the OPED algorithm:

\medskip\noindent
{\bf OPED Algorithm}.
{\it Let $m$ be a positive integer. For each reconstruction point $(x,y) \in B$, 
\begin{equation} \label{eq:A2m}
 \CA_{\rm OPED} : = \sum_{\nu=0}^{2m} \sum_{j=0}^{2m}
  g_{j,\nu} T_{j,\nu} (x,y), \qquad g_{j,\nu} =: 
         \CR f \left(\phi_\nu, \cos \psi_j \right), 
\end{equation}
where $\phi_\nu =\frac{2\nu \pi}{2m+1}$, $ \psi_j = 
\frac{(2j+1) \pi}{4m+2}$, and  
\begin{equation} \label{eq:Tjv}
T_{j,\nu} (x,y) =  \frac{ 1 }{(2m+1)^2} \sum_{k=0}^{2m} (k+1) 
   \sin(k+1)\psi_j U_k(x \cos \phi_\nu + y \sin \phi_\nu).
\end{equation}
}
\medskip

The above version of OPED with be refereed to as OPED of 
type I as it comes from Gaussian quadrature based on the 
zeros $\cos \frac{(2j+1)\pi}{4m+2}$, $0\le j \le 2m$,  of the 
Chebyshev polynomial of the first kind.  An alternative is the
OPED of type II, which uses Gaussian quadrature based on 
the zeros, $\cos \frac{j \pi}{2m+1}$, $1 \le j \le 2m$, of the 
Chebyshev polynomials of the second kind, for which the sum 
over $j$ in \eqref{eq:A2m} should start from $j =1$ and $\psi_j$ 
should be replaced by $\psi_j = \frac{j \pi}{2m+1}$.  These two 
choices lead to different scanning geometry of x-ray data, the 
advantage of the first type is explained in \cite{XTC}. Unless
specifically stated,  OPED in the following means type I. The
type II case will be mentioned whenever appropriate, the
proof  can be worked out in the same way and will be omitted. 

For each $(x,y)$ in the domain $B$, the algorithm produces a
number $\CA_{\rm OPED}$ which approximates the value of 
the image at that point. Typically a larger $m$ produces a better 
approximation. The values of $\CA_{\rm OPED}$ on the 
reconstruction points, say on a grid of pixel points, give the
reconstruction of the image. The convergence of the 
OPED algorithm is equivalent to the approximation property of 
the linear operator 
\begin{equation} \label{eq:Af}
\CA_{2m} f(x,y) = \sum_{\nu=0}^{2m} \sum_{j=0}^{2m} 
   \CR f(\phi_\nu,\cos\psi_j) T_{j,\nu}(x,y)
\end{equation}
where $f$ is the function that represents the image. The operator 
$\CA_{2m}$ preserves polynomials of degree $2m-1$, that is, 
$\CA_{2m} f \equiv f $ if $f$ is a polynomial of degree at most 
$2m-1$. In other words, the algorithm produces an image exactly 
if the image happens to be represented by a polynomial of degree 
less than $2m-1$. Furthermore, let  $\|\CA_{2m}\|$ denote the
operator norm of $\CA_{2m}$ in the uniform norm over $B$; 
then (\cite{X05})
\begin{equation}\label{eq:norm}
  \|\CA_{2m}\| = \CO(m \log (m+1) ), \qquad \hbox{as $m \to \infty$}.
\end{equation} 
As a consequence of this estimate, it follows that $\CA_{2m} f$ 
converges uniformly to $f$ on the disk $B$ if $f$ has continuous 
second order derivatives. The numerical test has shown that OPED 
gives fairly accurate reconstruction even when the data has sharp 
singularities. The proof of \eqref{eq:norm} is based on the following 
compact formula of $T_{j,\nu} (x,y)$. 

\begin{prop} \label{prop:Tjv}
Let $\theta_\nu: = \theta_\nu(x,y)$. Then
\begin{align*}
(2m+1)^2 T_{j,\nu}(x,y) = &
  \frac{(-1)^j (2m+1)\sin (2m+1)\theta_\nu}{2 \sin \theta_\nu}
 -  \frac{(-1)^j(2m+1)\cos (2m+1)\theta_\nu}
    {\cos \psi_j-\cos \theta_\nu} \\
 & - \frac{\sin \theta_\nu \sin \psi_j - (-1)^j\sin (2m+1)\theta_\nu
    (1- \cos\theta_\nu \cos\psi_j)}{2 \sin \theta_\nu
       (\cos \psi_j-\cos \theta_\nu)^2}. 
\end{align*}
\end{prop}

It should be mentioned that \eqref{eq:norm} was proved in 
\cite{X05} for OPED of type II. The compact formula of 
$T_{j,\nu}$ in \cite{X05} is also established for type II. So,
the formula of $T_{j,\nu}$ in the above proposition and 
\eqref{eq:norm} are in fact new, but the proof is very much 
similar to that of \cite{X05} and we choose not to repeat it. 

The use of orthogonal expansion in reconstructing functions 
from their Radon data can be traced back to the classical paper 
of Cormack, whose Radon inversion formula is based on the
expansions of $f$ and $Rf$ in spherical harmonics 
(\cite[p. 25]{N}).  The orthogonal polynomials on the disk have 
been used in the papers \cite{LS, Marr} and in the singular 
value decomposition of Radon transforms (see \cite{N}). The 
formula \eqref{S2m} and the resulted OPED algorithm are 
introduced in \cite{X05}. As we mentioned before that \eqref{S2m}
itself can be deduced from \cite{LS} as well as {BG,P}. The 
OPED of type II is closely related 
to an algorithm in \cite{BO}.  The connection to orthogonal 
polynomial expansion was not considered in \cite{BO}, nor was 
the convergence studied there. For other algorithms and image
reconstruction based on polynomial approximation, see 
\cite{LS, Ma, Marr, N} and the references therein. 

\subsection{Fast implementation of OPED algorithm}  

The structure of the $\CA_{2m}$ in \eqref{eq:A2m} allows us to use FFT 
(fast Fourier transform) once in a straightforward manner. In fact, 
let us define
$$
  S_{k,\nu} = \frac{k+1}{(2m+1)^2}\sum_{j=0}^{2m} g_{j,\nu} \sin(k+1)\psi_j.  
$$
Then $S_{k,\nu}$ can be evaluated by FFT for discrete sine transform. We 
can write $\CA_{\rm OPED}$ as 
\begin{equation}\label{eq:A2m2}
 \CA_{\rm OPED} = \sum_{\nu=0}^{2m} \sum_{k=0}^{2m}
   S_{k,\nu} U_k(x \cos \phi_\nu + y \sin \phi_\nu).
\end{equation}
Thus, the main step of the OPED algorithm lies in the evaluation of the 
above double sum, which can be considered as a back projection step. 
Let $N=2m+1$. Then the evaluation of the matrix $S_{k,\nu}$ costs 
$\CO(N^2\log N)$ operations (flops). The evaluation of the double sum costs 
$\CO(N^2)$ operations. Hence, the cost of evaluation on a grid of $M \times M$
is $\CO(N^2 (M^2 + \log N))$. In particular, if $M \approx N$, then the 
cost is $\CO(N^4)$. The main operation cost is at the evaluation of the 
double sums at the grid points.  In other words, the main cost lies in the 
back projection step. 

Unlike the FBP algorithm, the sum in \eqref{eq:A2m2} does not contain a 
discrete convolution that can be evaluated via FFT. However, the formula 
of  $U_k$ in \eqref{Cheby} allows us to write 
\begin{equation}\label{eq:A2m3}
  \CA_{\rm OPED}(x,y) = \sum_{\nu=0}^{2m} \frac{1}{\sin \theta_\nu}
    \sum_{k=0}^{2m} S_{k,\nu} \sin (k+1)\theta_{\nu}, 
\end{equation}
where 
$$
\theta_\nu : = \theta_\nu(x,y) = \arccos (x \cos \phi_\nu + y \sin \phi_\nu).
$$
The inner sum can be evaluated by FFT for discrete sine transform at certain 
points, which suggests that we introduce an interpolation step to take 
advantage of the fast evaluation by FFT. To be more precise, we define 
$$
\alpha_{\nu}(\theta)= \sum_{k=0}^{2m} S_{k,\nu} \sin (k+1) \theta, \qquad
   0 \le \theta \le \pi, 
$$
after the inner sum of \eqref{eq:A2m3}. The FFT for discrete sine transforms can 
be used to evaluate the numbers  
$$
\alpha_{l,\nu} : = \alpha_\nu( \xi_l), \qquad \xi_l := 
 \frac{(l+1) \pi} {2m+1}, \qquad l = 0 ,1, \ldots, 2m -1 
$$
effectively. That is, the inner sum in \eqref{eq:A2m3} can be evaluated 
by FFT when $\theta_\nu(x,y) = \xi_l$. For the interpolation step, we choose linear 
interpolation, which is also used in the implementation of FBP 
\cite[p. 109]{N}. For a given $(x,y)$, we choose the integer $l$ such 
that $\theta_\nu(x,y)$ lies between $\xi_l$ and $\xi_{l+1}$ and use 
the value of the linear interpolation between $\alpha_{l,\nu}$ and 
$\alpha_{l+1,\nu}$ as an approximation to the inner sum of \eqref{eq:A2m3}.
The linear interpolation is given by  
\begin{equation}\label{linear}
\ell_\nu(\theta) =  u_\nu(\theta)  \alpha_{l+1,\nu}
   + (1-u_\nu(\theta)) \alpha_{l,\nu}, \qquad 
      u_\nu(\theta): = \frac{ \theta - \xi_l}{ \xi_{l+1} - \xi_l},
\end{equation}
where $\xi_l \le \theta \le \xi_{l+1}$.  Then the fast implementation of OPED is given as follows: 

\medskip\noindent
{\bf Fast OPED algorithm}: Let $m$ be a positive integer.
 
\smallskip \noindent
 {\it Step 1.   For each $\nu = 0,\ldots, 2m$, use FFT to compute for each $k = 
0,\ldots, 2m$,
$$
 S_{k,\nu} = \frac{k+1}{(2m+1)^2}\sum_{j=0}^{2m} g_{j,\nu} \sin(k+1)\psi_j  .
$$
Step 2. For each $l=0,1,\ldots, 2m-1$, use FFT to compute 
$$
\alpha_{l,\nu}: = \sum_{k=0}^{2m} S_{k,\nu} \sin (k+1)\xi_l.
$$
Step 3. For each reconstruction point $(x,y)$ inside the disk of the
radius $\cos \frac{\pi}{2m+1}$, determine integers $l$ such that
$$
 l = \left \lfloor\frac{2m+1}{\pi} \theta_\nu \right \rfloor -1,
   \qquad  \hbox{where} \quad
  \theta_\nu = \arccos (x \cos \phi_\nu + \sin \phi_\mu),
$$
and evaluate  
$$
 f_{\rm OPED} = \sum_{\nu = 0}^{2m} \frac{1}{\sin \theta_\nu} \left[ 
      (1- u_\nu ) \alpha_{l,\nu} +  u_\nu \alpha_{l+1,\nu} \right],        
$$
where $u_\nu = (2m+1) \theta_\nu / \pi - (l + 1)$. 
 }
\medskip

A fast implementation for OPED of type II works similarly, with 
$\xi_l$ replaced by $\xi = \frac{(l +\frac{1}{2})\pi}{2m+1}$, $0\le 
l \le 2m$. In other words, exchanging the values of $\psi_j$ 
and $\xi_l$ in the above fast algorithm for OPED of type I leads to
the fast algorithm for OPED of type II. 

A couple of remarks are in order. First of all, $\sin \theta_\nu$ 
appears in the denominator in the last step of the algorithm. 
However, $\sin \theta_\nu  =0$ only if $\cos \theta_\nu (x,y) = 
x \cos \phi_\nu + y \sin \phi_\nu =1$, which happens only if 
$(x,y) = (\cos \phi_\nu, \sin\phi_\nu)$. Since the points 
$(\cos \phi_\nu, \sin\phi_\nu)$ are on the boundary of the 
region $B$, we do not have to take them as reconstruction 
points.  In fact, the region of interests is usually inside the unit 
disk; thus, we can evaluate at points inside a smaller disk in 
$B$, and this will also ensure that the values of $\sin \theta_\nu$ 
in the last step will not be too small to cause loss of significant 
digits in the computation. Furthermore, if we restrict $(x,y)$ to 
a disk with radius $\cos \xi_0 = \cos \pi/(2m+1)$, then it will also 
guarantee that the choice of $l$ in the Step 3 is unique for 
all $(x,y)$ in that disk. 

The algorithm uses FFT twice, the final sum in Step 3 is a single 
sum whose evaluation costs $\CO(N)$ operations. Hence, the 
cost of evaluations on an $M\times M$ grid with $M \approx N$ 
is $\CO(N^3)$, which is in the same order of magnitude as the 
FBP algorithm. In fact, the structure of the OPED algorithm with 
this fast implementation is similar to that of FBP algorithm 
(cf. \cite[p. 109]{N}). 

\section{Convergence of OPED algorithm with linear Interpolation}
\setcounter{equation}{0}

Throughout this section, we let $\|\cdot\|$ denote the uniform norm 
on $B$.

As we mentioned before that the convergence of the OPED algorithm depends on the norm $\CA_{2m}$ of the operator defined in 
\eqref{eq:Af}. Following the proof in \cite{X05}, we have 
$$
  \|\CA_{2m}\| = \max_{(x,y) \in B} \sum_{\nu = 0}^{2m} \sum_{j=0}^{2m}
      \sin \psi_j  | T_{j,\nu}(x,y)| = \CO(m \log (m+1)).
$$ 
As a consequence of this estimate and the fact that $\CA_{2m}$ 
preserves polynomials of degree $2m-1$, the triangle inequality implies 
that 
\begin{equation}\label{conv}
 \|\CA_{2m}f - f\| \le c \, m \log (m+1) E_{2m-1}(f) 
\end{equation}
where $E_n(f) = \inf \{\|f - P_n|: \deg P_n \le n\}$ is the error of best 
approximation to $f$ by polynomials of degree at most $n$. It is known
(\cite{X05a}) that if $f \in C^{2r}(B)$, then  
\begin{equation} \label{best}
     E_n(f) \le c n^{- 2r} \|\CD^r f\|, 
\end{equation}
where $\CD$ is a differential operator of order $2r$.  In particular,
\eqref{conv} and \eqref{best} show that $\CA_{2m} f$ converges uniformly 
to $f$ if $f$ has continuous second order derivatives. 

When the interpolation step is used in  the evaluation of OPED, 
the operator is changed and we denote the new operator by 
$\CA\CI_{2m}$, which is given by 
$$
\CA\CI_{2m} f(x,y):=\sum_{\nu = 0}^{2m} \frac{1}{\sin (\theta_\nu(x,y))}
     \ell_\nu (\theta_\nu(x,y)),
$$
where $\ell_\nu$ is defined in \eqref{linear}. It turns out that the 
norm of the operator $\CA\CI_{2m}$ has the same growth order 
as that of $\CA_{2m}$. 

\begin{prop} \label{prop:3.1}
For $m \ge 0$, let $\Omega_m =\{(x,y): \sqrt{x^2+y^2} \le 
\cos \pi/(2m+1)\}$. Then 
$$
 \max_{ (x,y) \in \Omega_m} | \CA\CI_{2m} f(x,y) | \le c \|f\| 
    \, m (\log (m+1)).
$$
\end{prop}

\begin{proof}
Let $u_\nu = (\theta_\nu - \xi_l)/(\xi_{l+1} - \xi_l)$. The way that the 
index $l$ is chosen implies that $0 \le u_\nu \le 1$. Recall the 
definition of $T_{j,\nu} (x,y)$ in \eqref{eq:Tjv}. We shall abuse the 
notation somewhat and write $T_{j,\nu} (\xi_l)$ when 
$\arccos(x\cos \phi_\nu+ y \sin \phi_\nu)=\xi_l$. 
Using the definition of $U_n(t)$ in \eqref{Cheby}, the operator 
can be written as
\begin{align*}
\CA\CI_{2m} f(x,y) = & \sum_{\nu=0}^{2m} \frac{1}{\sin \theta_\nu(x,y)}
  \sum_{j=0}^{2m} \CR f(\phi_\nu,\cos\psi_j) \\
   & \quad  \times \left[(1-u_\nu) T_{j,\nu}(\xi_l) \sin \xi_l + 
      u_\nu T_{j,\nu}(\xi_{l+1}) \sin \xi_{l+1} \right].
\end{align*}
For $(x,y) \in \Omega_m$, write $x = r \cos \phi$ and 
$y = r \sin \phi$, then $x\cos\phi_\nu + y\sin \phi_\nu =  
 r \cos (\phi-\phi_\nu) \le r$ and $r \le \cos \pi/(2m+1)$. In particular,  
\begin{equation} \label{theta}
 \sin \theta_\nu(x,y) = \sqrt{1- \cos^2 \theta_\nu(x,y)} \ge
   \sqrt{1-\cos^2 \tfrac{\pi}{2m+1}}  = \sin \tfrac{\pi}{2m+1}.
\end{equation}
Furthermore, the definition of $\xi_l$ shows that there is a constant 
$c$ independent of $m$ such that 
$$
\left|\frac{\sin \xi_l}{\sin (\theta_\nu(x,y))} \right| \le c  
\qquad \hbox{and} \qquad  
\left|\frac{\sin \xi_{l+1}}{\sin (\theta_\nu(x,y))} \right| \le c 
$$
for $(x,y) \in \Omega_m$. Consequently, we conclude that 
\begin{equation} \label{CAI}
 \left|\CA\CI_{2m} f(x,y)\right| \le 
 \sum_{\nu=0}^{2m} \sum_{j=0}^{2m} |\CR f(\phi_\nu,\cos\psi_j)|
    \left( |T_{j,\nu}(\xi_l)| +| T_{j,\nu}(\xi_{l+1})| \right).
\end{equation}
Using the fact that $|\CR f(\phi,t)| \le \sqrt{1-t^2} \|f\|$ (see, 
for example, \cite{X05}), it follows that
$$
   \left|\CA\CI_{2m} f(x,y)\right| \le \|f\| 
 \sum_{\nu=0}^{2m} \sum_{j=0}^{2m} \sin \psi_j
    \left( |T_{j,\nu}(\xi_l)| +| T_{j,\nu}(\xi_{l+1})| \right).
$$
We now use the compact formula of $T_{j,\nu}$ in 
Proposition \ref{prop:Tjv}. Since $\sin (2m+1)\xi_l = 0$ and 
$\cos (2m+1) \xi_l = (-1)^{l+1}$, it follows that 
$$
(2m+1)^2 T_{j,\nu}(\xi_l) = 
   - \frac{(-1)^{j+l+1}(2m+1)}{\cos \psi_j-\cos \xi_l} 
    - \frac{(-1)^{l+1} \sin \psi_j}{2(\cos \psi_j-\cos \xi_l)^2}. 
$$
Since $\phi_j = \frac{(j+1/2)\pi}{2m+1}$ and $\xi_l = 
\frac{(l+1)\pi}{2m+1}$,  the denominator of $T_{j,\nu}(\xi_l)$ is 
never zero. We have 
$$
   \cos \psi_j - \cos \xi_l = 2 \sin \frac{\xi_l - \psi_j}{2}
      \sin \frac{\xi_l+ \psi_j}{2}.
$$
Since $\psi_{2m-j} = \pi - \psi_j$ and $\xi_{2m-l} = \pi - \xi_{l-1}$,
we can assume that $0 < \xi_l < \pi/2$, which means that 
$0 \le l \le m -1$.  If $0\le \psi_j \le \pi/2$, we have 
$$
  \frac{ \sin \psi_j}{\sin \frac{\psi_j+\xi_l}{2}} = 2
      \frac{ \sin \frac{\psi_j}{2} \cos \frac{\psi_j}{2}}
          {\sin \frac{\psi_j+\xi_l}{2}} \le  2 \cos \frac{\psi_j}{2} \le 2. 
$$
The above equation also holds if $\pi/2 < \psi_j < \pi$, since then 
$\pi/4 \le \frac{\psi_j + \xi_l}{2} \le 3 \pi/4$ and 
$\sin \frac{\psi_j + \xi_l}{2}  \ge \sqrt{2}/2$. Hence, using the fact that
$\sin \theta \ge (2/\pi) \theta $ for $0\le \theta \le \pi/2$ we have 
\begin{align*}
\sum_{\nu = 0}^{2m} \sum_{j=0}^{2m} \sin \psi_j |T_{j,\nu}(\xi_l)|
 & \le \sum_{j=0}^{2m} \frac{1}{|\sin \frac{\psi_j - \xi_l}{2}|}
  + \frac{1}{2m+1} \sum_{j=0}^{2m}\frac{1}{|\sin \frac{\psi_j - \xi_l}{2}|^2}\\
 & \le \pi \sum_{j=0}^{2m} \frac{1}{|\psi_j - \xi_l|} 
      + \frac{\pi^2}{2m+1} \sum_{j=0}^{2m}\frac{1}{|\psi_j - \xi_l|^2} \\
 & \le \sum_{j=0}^{2m} \frac{2m+1}{|j - l + 1/2|} 
      + \sum_{j=0}^{2m}\frac{2m+1}{|j - l + 1/2|^2} \\
 & \le c (2m+1) \log (m+1).
\end{align*}
The sum involving $T_{j,\nu}(\xi_{l+1})$ is estimated in exactly the same way.
By \eqref{CAI}, the proof is completed. 
\end{proof} 

This proposition shows that the interpolation step does not 
increase the growth order of the operator norm, at least when 
we restrict the norm to the region $\Omega_m$. Unlike 
$\CA_{2m}$, the operator $\CA\CI_{2m}$ no longer preserves 
polynomials of high degrees. We need an estimate of the 
error $\CA\CI_{2m} f -f $ for 
$f$ being a polynomial. First, however, we need some results on 
$\proj_k: L^2(B) \mapsto \CV_k(B)$ defined in \eqref{OPexpan}. Let 
$P_n(\xb,\yb)$, $\xb,\yb \in B$,  denote the reproducing kernel of $\CV_k(B)$. 
Then 
$$
   \proj_k f(\xb) =  \frac{1}{\pi} \int_B f(\yb) P_k(\xb,\yb) d\yb. 
$$
Furthermore, the kernel $P_k$ satisfies a compact formula (\cite{X99})
\begin{equation}\label{reprod}
  P_k(\xb,\yb) =  \frac{k+1}{2} \int_{0}^\pi U_k \left(\langle \xb, \yb\rangle + 
     s \sqrt{1-\|\xb\|^2}\sqrt{1-\|\yb\|^2} \right) \frac{ds}{\sqrt{1-s^2}},
\end{equation}
where $\langle \cdot,\cdot \rangle$ and $\|\cdot\|$ are the usual Euclidean
inner product and norm on $\RR^2$, and $U_k$ is the Chebyshev 
polynomial of the second kind.  

\begin{lem} \label{lem1}
Let $\|f\|_2$ denote the $L^2(B)$ norm of $f$. Then the projection operator 
satisfies 
$$
    |\proj_kf(\xb)| \le   (k+1) \|f\|_2, \qquad \forall \xb \in B.
$$
\end{lem}

\begin{proof}
Using the Cauchy-Schwartz inequality and the reproducing property of 
$P_k(\cdot,\cdot)$, we have
$$
    |\proj_k f(\xb)|  \le \|f\|_2 \left(\frac{1}{\pi} \int_{B} |P_k(\xb,\yb)|^2 
       d\yb \right)^{1/2}  =  \|f\|_2 \left [ P_k(\xb,\xb) \right ]^{1/2}. 
$$
The definition \eqref{Cheby} gives the well known inequality $|U_k(t)| \le k+1$,
so that by \eqref{reprod} we have $|P_k(\xb,\xb)| \le (k+1)^2$, from which 
the stated result follows. 
\end{proof}

Using $\proj_k f$ we can construct a sequence of 
polynomials that approximates $f$. For this we let $\eta$ be a nonnegative 
$C^\infty$ function on $\RR$ that satisfies
$$
\eta(t) = 1, \quad 0 \le t \le 1, \qquad \hbox{and} \quad \supp \eta (t) = [0,2].
$$
For $f \in L^2(B)$ we define $S_n^\eta f$ by 
$$
   S_n^\eta f(\xb) = \sum_{j=0}^{2n} \eta\left(\frac{j}{n}\right) \proj_jf(\xb).
$$
Let $\Pi_n^2$ denote the space of polynomials of degree at most $n$ in two
variables. The approximation property of $S_n^\eta$ is given in the following
lemma (\cite{X05a}). 

\begin{lem} \label{lem2}
Let $f \in C(B)$. Then
\begin{enumerate}
\item $S_n^\eta f \in \Pi_{2n-1}^2$ and $S_n^\eta P  = P$ for $P\in \Pi_n^2$;
\item for $n \in \NN$,
$$
    \|S_n^\eta f \| \le c \|f\| \quad \hbox{and}\quad  \|f- S_n^\eta f \| \le c E_n (f). 
$$
\end{enumerate}
\end{lem} 

This shows that $S_n^\eta f$ realizes, up to a constant, the best approximation 
by polynomials. For this sequence of polynomials, we have the following: 

\begin{lem} \label{prop:3.2}
Let $F = S_q^\eta$ and $q < m$. For each compact subset $\Omega$ in $B$,
there is a constant $c$, independent  of $q$, $m$ and $f$, such that
$$
 \max_{(x,y) \in \Omega}  |\CA\CI_{2m}F (x,y) - F (x,y)| 
           \le  c\,\|f\| \frac{q^4}{ (2m+1)^2}.
$$
\end{lem}

\begin{proof}
Since $F$ is a polynomial of degree $2q$, we write its expansion as
$$
    F(x,y)  = \sum_{j=0}^{2q}  F_j(x,y), \qquad F_j = \eta\left(\frac{j}{q}\right)
           \proj_j f(x) \in \CV_j(B). 
$$
The Radon projections of orthogonal polynomials can be explicitly computed.
We have  (\cite{Marr}),
$$
\CR F_j (\phi,t) = \frac{2}{j+1} \sqrt{1-t^2} U_j(t) F_j(\cos \phi, \sin\phi). 
$$
The polynomials $U_k(t)$ are orthogonal with respect to $\sqrt{1-t^2}$ on
the interval $[-1,1]$. Consequently, 
$$
 \frac{1}{\pi} \int_{-1}^1 \CR F_j (\phi, t) U_k(t) dt 
       = \frac{1}{k+1}  F_j (\cos \phi,\sin \phi) \delta_{k,j}. 
$$
Hence, using \eqref{S2m} and the fact that $q< m$,  we conclude that 
\begin{align} \label{Af}
 F (x,y) =  & S_{2m}F(x,y)  = \sum_{j=0}^{2q} S_{2m} F_j (x,y)\\ 
  = & \frac{1}{2m+1}\sum_{\nu=0}^{2m} \sum_{j=0}^{2q} F_j(\cos \phi_\nu,\sin\phi_\nu)
              U_j(\cos \theta_\nu(x,y)). \notag
\end{align}
In exactly the same way, we also obtain that 
\begin{align} \label{AIf}
 \CA\CI_{2m} f(x,y) =   \frac{1}{2m+1}\sum_{\nu=0}^{2m} & \frac{1}{\sin \theta_\nu(x,y)}
   \sum_{j=0}^{2q} F_j(\cos \phi_\nu,\sin\phi_\nu) \\
   & \times         \left[(1-u_\nu)\sin (j+1)\xi_{l+1}+ u_\nu\sin(j+1)\xi_l\right], \notag
\end{align}
where $u_\nu = u_\nu(\theta_\nu(x,y))$ is the linear function in \eqref{linear}. 
Introducing the function
$$ 
      H_\nu(\theta) = \sum_{j=0}^{2q}  F_j(\cos \phi_\nu,\sin\phi_\nu) \sin(j+1)\theta,
$$
which depends on $F$, we then conclude that 
\begin{align*}
  F(x,y) & - \CA\CI_{2m} F(x,y) \\
     &= \frac{1}{2m+1} \sum_{\nu=0}^{2m} \frac{1}{\sin \theta_\nu(x,y)}
       \left[ H_{\nu}(\theta_\nu(x,y)) - (1-u_\nu) H_{\nu} (\xi_l) - 
              u_\nu H_{\nu} (\xi_{l+1})\right].
\end{align*}
The expression in the square bracket is the difference of $H_\nu(\cdot)$ 
and its linear interpolation at $\xi_l$ and $\xi_{l+1}$ evaluated at 
$\theta_\nu(x,y)$.
Hence the well-known estimate of the error of linear interpolation shows 
\begin{align*}
   &\left[ H_{\nu}(\theta_\nu(x,y)) - (1-u_\nu) H_{\nu} (\xi_l) - 
              u_\nu H_{\nu} (\xi_{l+1})\right] \\
    & \le \frac{1}{2} \|H_\nu''\| |(\theta_\nu(x,y) -\xi_l)(\theta_\nu(x,y) -\xi_{l+1})| 
    \le \frac{1}{16} \|H_\nu''\| \frac{\pi^2}{(2m+1)^2},
\end{align*}
where $\|H_\nu ''\|= \max_{0\le \theta \le \pi}| H_\nu''(\theta)|$. Since $H_\nu$ 
is a trigonometric polynomial of degree $2q$, the classical Bernstein inequality 
(cf. \cite[Vol. 2, p. 11]{Z}) shows that 
$$
     \|H_\nu''\| \le (2q)^2 \|H_{\nu}\| := 2q^2 \max_{0\le t \le \pi} | H_\nu(t)|.
$$
Furthermore, by the definition of $H_\nu$ and Lemma \ref{lem1}, we have 
\begin{align*}
    \|H_{\nu} \| \le \sum_{j=0}^{2q} |F_j(x,y)| & \le \sum_{j=0}^{2q} |\proj_j F(x,y)| \\
           &  \le  \|f\|_2 \sum_{j=0}^{2q} (j+1) \le  2(q+1)^2 \|f\|. 
\end{align*}
Hence, putting these inequalities together, we conclude that 
$$
   |F(x,y) - \CA\CI_{2m} F(x,y)| \le c  \|f\| \frac{q^4}{(2m+1)^3} \sum_{\nu=0}^{2m}  
            \frac{1}{\sin \theta_\nu(x,y)} . 
$$ 
If $(x,y) \in \Omega$, a compact set in $B$, then there is an $r < 1$ such
that $\sqrt{x^2 + y^2 } \le r$, so that $\sin\theta_{\nu}(x,y) \ge \sqrt{1-r^2}$ 
as in \eqref{theta}. Thus the last
sum is bounded by a constant multiple of $m$ and the stated result follows.
\end{proof} 

\begin{thm}
If $f \in C^4(B)$ and $\Omega$ is a compact subset of $B$, then there is 
a constant $c_f$, independent of $m$, such that 
\begin{equation} \label{bound}
 \max_{(x,y) \in \Omega}  | \CA\CI_{2m} f(x,y) - f (x,y)| \le c_f \frac{\log (m+1)}{\sqrt{m}}.
\end{equation} 
 In particular, $\CA\CI_{2m} f$ converges uniformly to $f$ on any compact subset 
 of $B$. 
\end{thm}

\begin{proof}
Let $F = S_q^\eta f$ and $q < m$ as in the previous lemma and 
$\Omega$ be a compact subset of $B$. Since $\CA\CI_{2m}$ is a 
linear operator, by the 
triangle inequality, Proposition \ref{prop:3.1}, Lemma \ref{lem2} and 
Lemma \ref{prop:3.2},  we obtain
\begin{align*}
  |\CA\CI_{2m} f(x,y) - f(x,y)| & \le (1+ \|\CA\CI_{2m}\|) |(f-F)(x,y)| +
         |\CA\CI_{2m}F(x,y) - F(x,y)| \\
  &  \le  c \left( m \log (m+1) E_q(f) +   m^{-2} q^4 \|f\| \right)
\end{align*}
for all $(x,y) \in \Omega$ and $q < m$. In particular,  if $f \in C^4(B)$, then
 \eqref{best} implies that 
\begin{equation} \label{estimate}
 \max_{(x,y)\in \Omega} |\CA\CI_{2m} f(x,y) - f(x,y)| 
    \le  c \left( m \log (m+1) q^{-4} \|\CD^2 f\| +   m^{-2} q^4 \|f\| \right).
\end{equation}
In particular, setting $q \approx m^{3/8}$ in the above inequality gives 
\eqref{bound}. 
\end{proof}

The theorem is stated for functions in $C^4(B)$, which is likely 
too restrictive and the convergence could hold under less 
restrictive conditions.  Also, for functions that are more smooth, 
one could get better convergence rate by choosing $q$ in 
\eqref{estimate} differently. The rate so obtained and the one 
in \eqref{bound}, however, are likely not sharp and reflect the 
limitation of our method of proof.  

The main merit of the theorem is that it establishes the 
convergence of OPED algorithm with linear interpolation for 
smooth functions (images).  In practice, the images often have 
sharp edges, which means that the functions representing 
images may not be smooth or even continuous.  In the following 
section we present numerical examples, which demonstrate 
that $\CA\CI_{2m}f $, OPED with interpolation step, converges 
well even for functions that are not continuous and it converges 
almost as well as $\CA_{2m} f$, OPED without interpolation step. 

\section{Implementation and Result}
\setcounter{equation}{0}

For numerical implementation, we used the FFT for discrete sine 
transform in the package FFTW (http://www.fftw.org/). The
numerical example is conducted on the Shepp-Logan head 
phantom \cite{SL} (see Figure 2).  This is an analytic phantom, 
highly singular, as the image contains jumps at the boundary of 
every ellipse in the image, include the one on the boundary. The 
function that represents  the image is a step function, which is not 
continuous at the boundary of each of the ellipses. 

We reconstruct the image with  OPED algorithm without the 
interpolation step and Fast OPED algorithm, which contains 
the interpolation step as shown in previous section, respectively. 
In both cases,  $S_{k,\nu}$ are computed with FFT. 

In Figure 1 images reconstructed by the original OPED and Fast 
OPED algorithms are depicted side by side. In these images we 
choose $m = 512$ and the size of the images are 
$512 \times 512$ pixels.  

\medskip

\begin{center}
\includegraphics[width = 6cm]{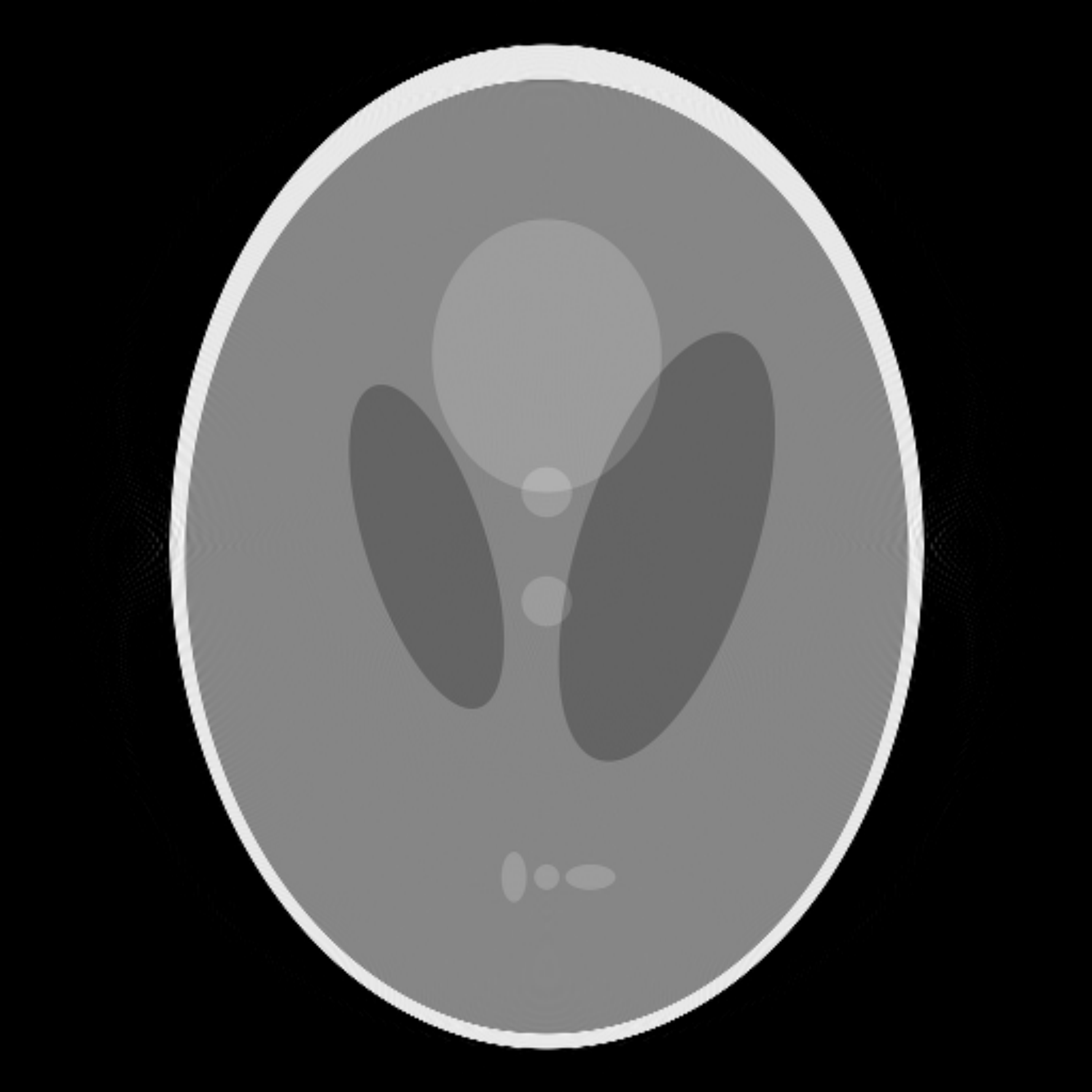} \hskip .2in
\includegraphics[width = 6cm]{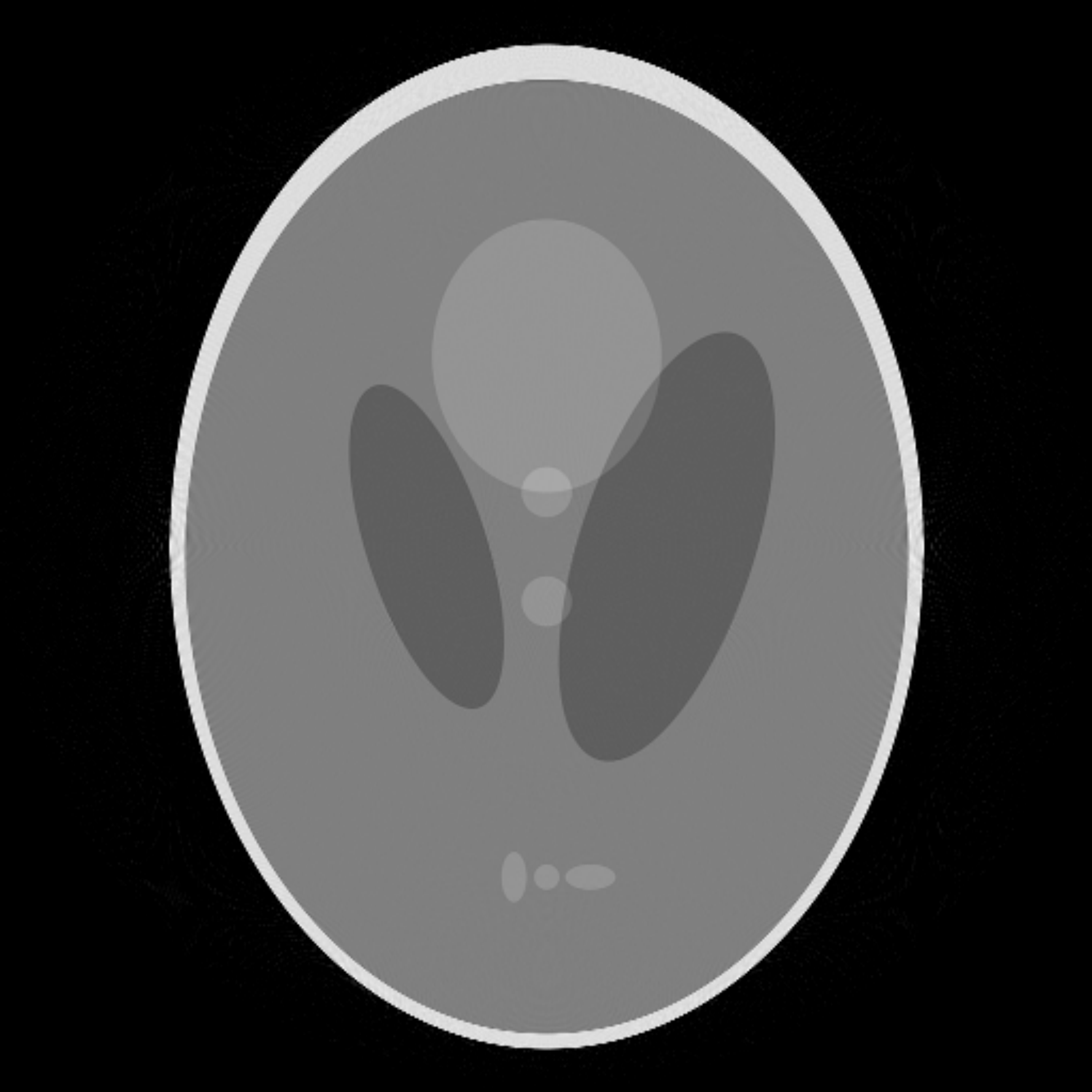} 
\end{center}
\noindent
Figure 1. Left: reconstruction by OPED algorithm with $m=512$. 
Right: reconstruction by Fast OPED algorithm with $m=512$. 

\medskip

The reconstruction is carried out on a CELSIUS R610 computer with two 
Intel Xeon(TM) CPU, each 3065 MHz, and 4 GB RAM. The code is written 
in C language. Using OPED algorithm,   the reconstruction took 344 seconds, 
in which more than 95\% of the time is used on the back projection step. 
Using Fast OPED algorithm, the reconstruction took merely 13 seconds, 
an improvement of more that 26 times.  Furthermore, the two images show 
almost no visual difference.  In Figure 2, the original Shepp-Logan phantom 
and the difference of the two images in Figure 1 are depicted. 
 
\medskip

\centerline{
\includegraphics[width = 6cm]{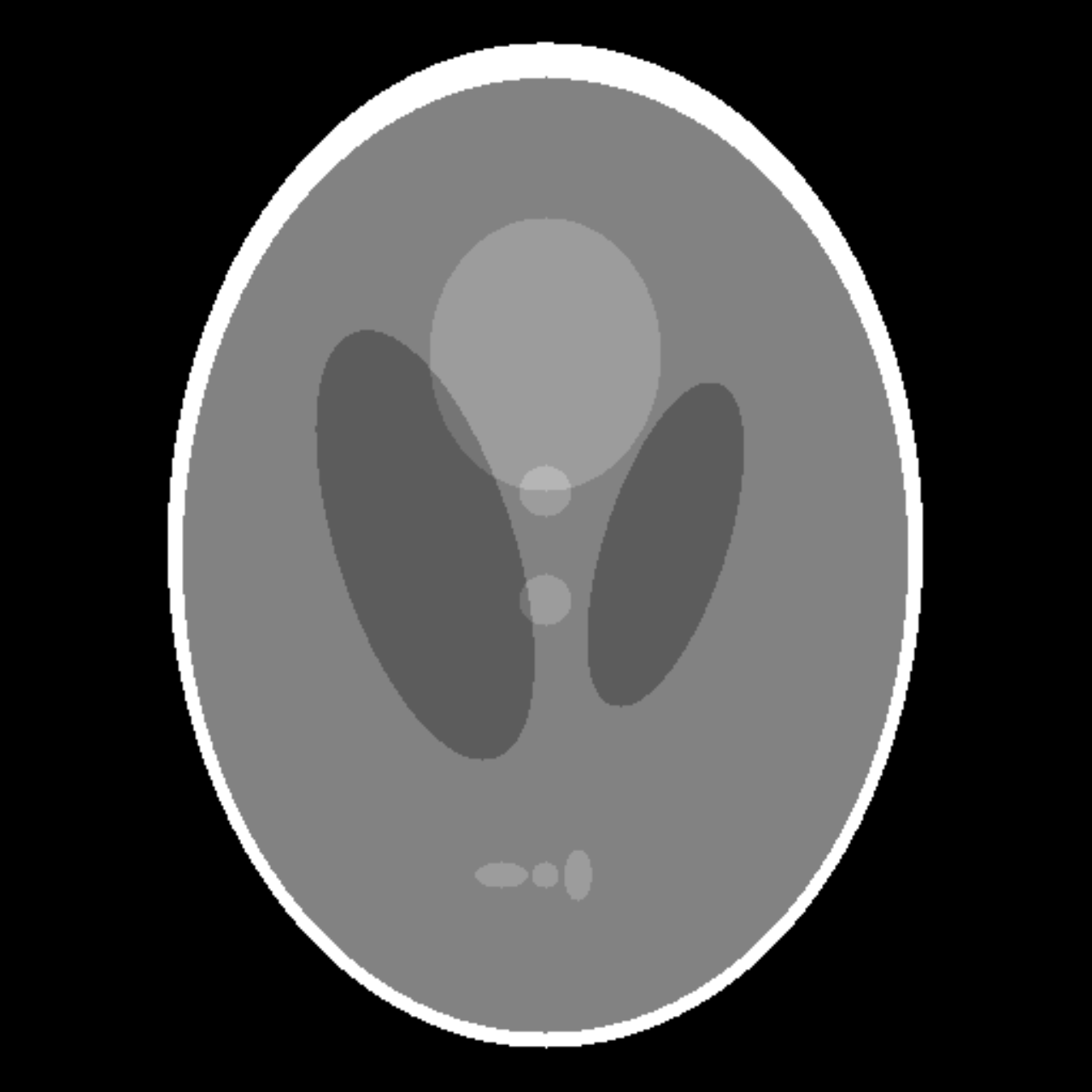} \hskip .2in
\includegraphics[width = 6cm]{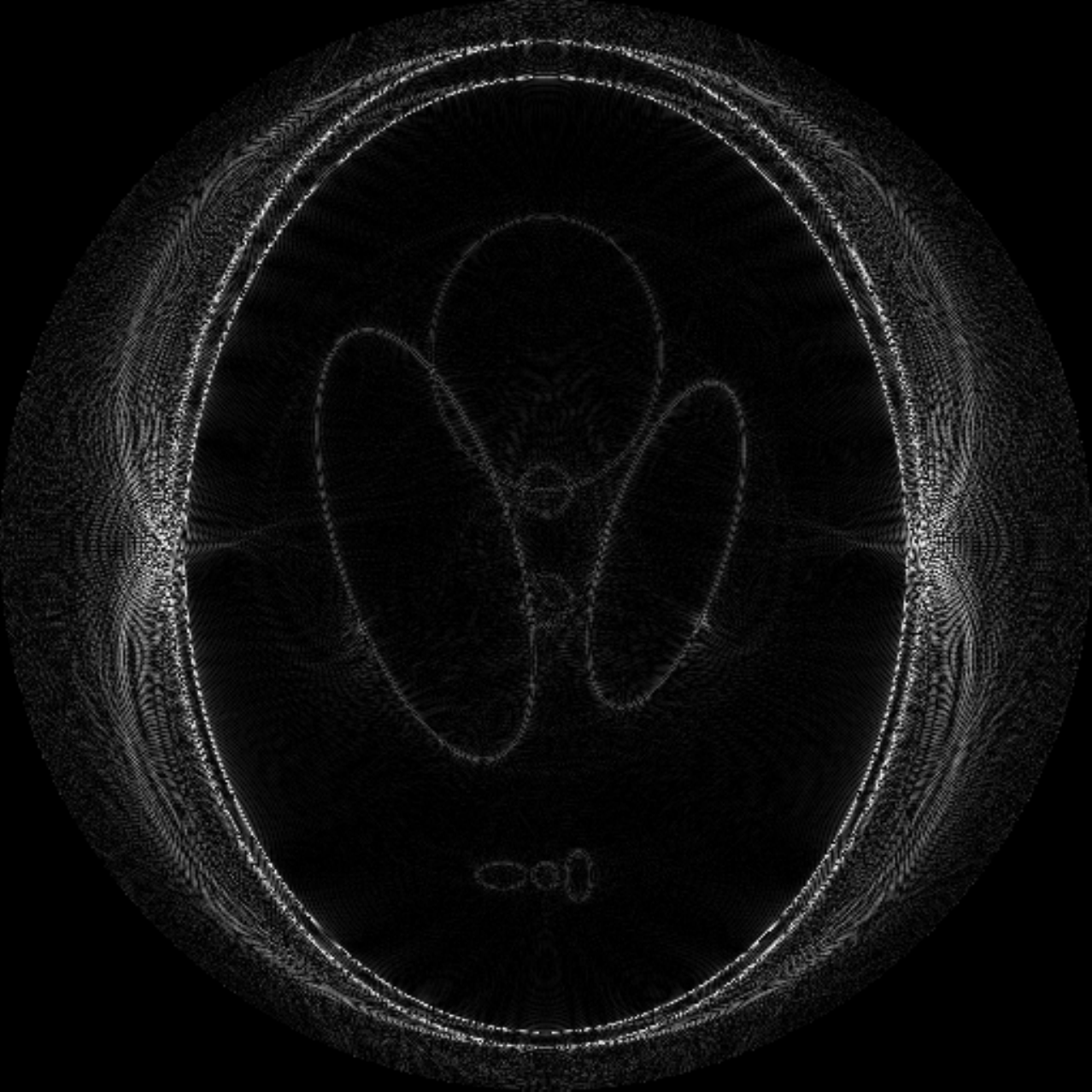}
}
\noindent
Figure 2. Left: original Shepp-Logan phantom. Right: the difference between the two reconstructed images in Figure 1. 

\medskip
 
We can also measure the errors of the reconstruction. Let $X$ denote an 
image, represented by its pixel values, so that  $X = \{X_i: 1 \le i \le N\}$, 
where $N$ is the number of pixels in the image.  Let $X^R$ denote the
reconstructed  image, $X^R = \{X_i^R: 1 \le i \le N\}$. The relative square 
total error (RSE) between $X$ and $X^R$ is defined by
$$
  Q_{RSE}= \frac{\sum_i (X_i^R - X_i)^2} {\sum_{i}^N (X_i^R)^2}, 
$$
and the mean error (ME) between  $X$ and $X^R$ is defined by 
$$
  Q_{ME} =  \frac{1}{N} \sum_{i=1}^N | X_i - X_i^R |.
$$
In our case $N = 512 \times 512 = 512^2$.  The result is reported in Table 1, 
in which {\sc orig} stands for  the original Shepp-Logan phantom, OPED and 
Fast OPED stand for reconstructed image by OPED and by Fast OPED, 
respectively.  For example, {\sc orig} vs OPED means the error between 
the original phantom and the reconstruction by OPED. 

\medskip
\begin{table}[h]
\begin{center} 
\begin{tabular*} {0.91\textwidth}
{   l  |c |c |c| r }
& {\sc orig} vs OPED     &  \sc{ orig}  vs {\sc Fast} OPED  
& OPED vs {\sc Fast} OPED \\
   \hline  
  RSE  & 0.00239702 & 0.00249574 & 0.000515499 \\  
    \hline
  ME   & 0.0129175  & 0.00981329 & 0.007715128   
\end{tabular*}
\caption{Error Estimates of OPED and Fast OPED}
\end{center}
\end{table}

\noindent
The results in the table shows that fast OPED is slightly worse in
relative least square error, but slightly better in mean error. The 
order of magnitude of the error is the same. The difference is practically 
negligible.


\section{Conclusion}
\setcounter{equation}{0}

We introduced a fast implementation of OPED algorithm by 
using FFT and a linear interpolation step. The fast algorithm is 
proved to converge uniformly on any compact subset in the disk 
and the numerical test has shown that it reconstructs images 
accurately and is as good as the original OPED algorithm.

In conclusion, the fast OPED algorithm with interpolation step 
is not only much faster, it also shares main merits of the original 
OPED algorithm. Hence, for practical applications, the fast 
OPED algorithm should be the recommended method.

\end{document}